\documentclass{amsproc}
\usepackage{amsmath}[1996/11/01]
\usepackage{amssymb,amsthm,amsxtra}

\def\[#1\]{\begin{equation}#1\end{equation}}
\makeatletter
\def\beq{%
   \relax\ifmmode
      \@badmath
   \else
      \ifvmode
         \nointerlineskip
         \makebox[.6\linewidth]%
      \fi
      $$
   \fi
}
\def\eeq{%
   \relax\ifmmode
      \ifinner
         \@badmath
      \else
         $$
      \fi
   \else
      \@badmath
   \fi
   \ignorespaces
}

\def\enddisplaymath{\eeq\global\@ignoretrue}
\makeatother

\newtheorem{thm}{Theorem}
\newtheorem{cor}[thm]{Corollary}
\newtheorem{lem}[thm]{Lemma}

\theoremstyle{remark}
\newtheorem*{rem}{Remark}
\newtheorem{rems}[thm]{Remark}

\theoremstyle{definition}
\newtheorem{defn}{Definition}

\numberwithin{equation}{section}
\numberwithin{thm}{section}
\numberwithin{eg}{section}

\newcommand{\la}{\langle}
\newcommand{\ra}{\rangle}

\begin{document}

\title{A difference-integral representation of Koornwinder polynomials}
\author{Eric M. Rains}
\address{Department of Mathematics, University of California, Davis}
\email{rains@math.ucdavis.edu}

\date{March 18, 2004}
\copyrightinfo{2004}{Eric M. Rains}

\begin{abstract}
We construct new families of ($q$-) difference and (contour) integral
operators having nice actions on Koornwinder's multivariate orthogonal
polynomials.  We further show that the Koornwinder polynomials can be
constructed by suitable sequences of these operators applied to the
constant polynomial $1$, giving the difference-integral representation of
the title.  Macdonald's conjectures (as proved by van Diejen and Sahi) for
the principal specialization and norm follow immediately, as does a
Cauchy-type identity of Mimachi.
\end{abstract}
\maketitle
\tableofcontents

\section{Introduction}

In \cite{KoornwinderTH:1992}, Koornwinder introduced a family of (symmetric)
multivariate orthogonal (Laurent) polynomials orthogonal with respect to
the following density on the unit torus:
\begin{align}
\Delta^{(n)}&(z_1,z_2,\dots z_n;t_0,t_1,t_2,t_3;q,t)\\
&=
\prod_{1\le i\le n}
\frac{(z_i^{\pm 2};q)}{(t_0 z_i^{\pm 1},t_1 z_i^{\pm 1},t_2 z_i^{\pm 1},t_3
  z_i^{\pm 1};q)}
\prod_{1\le i<j\le n}
\frac{(z_i^{\pm 1} z_j^{\pm 1};q)}
{(t z_i^{\pm 1} z_j^{\pm 1};q)},\notag
\end{align}
where $(x,y,z,\dots,w;q)$ represents the infinite $q$-symbol
\begin{align}
(x;q) &:= \prod_{j\ge 0} (1-q^j x),\\
(x,y,z,\dots,w;q)&:=(x;q)(y;q)(z;q)\cdots(w;q),
\end{align}
so in particular $(z_i^{\pm 1}z_j^{\pm
  1};q)=(z_iz_j;q)(z_i/z_j;q)(z_j/z_i;q)(1/z_iz_j;q)$.

To be precise, the Koornwinder polynomials $K^{(n)}_\lambda(\dots
z_i\dots;t_0,t_1,t_2,t_3;q,t)$ are uniquely defined by the following
requirements:
\begin{enumerate}
\item $K^{(n)}_\lambda(;t_0,t_1,t_2,t_3;q,t)$ is a $BC_n$-symmetric
  polynomial; i.e., a Laurent
polynomial invariant under permutations of the variables and substitutions
$z_i\mapsto z_i^{-1}$.
\item Moreover, it is monic with respect to dominance:
\beq
K^{(n)}_\lambda(\dots z_i\dots;t_0,t_1,t_2,t_3;q,t)
=
m_\lambda + \text{dominated terms}.
\eeq
\item With respect to the above density, it is orthogonal to any strictly
  dominated monomial.
\end{enumerate}
When $n=1$, Koornwinder's density becomes the following density associated
to the Askey-Wilson polynomials \cite{AskeyR/WilsonJ:1985}:
\[
\Delta^{(1)}(z;t_0,t_1,t_2,t_3;q,t)
=
\frac{(z^{\pm 2};q)}{(t_0 z^{\pm 1},t_1 z^{\pm 1},t_2 z^{\pm 1},t_3
  z^{\pm 1};q)}
\]
and thus the Koorwinder polynomials are a multivariate analogue of
Askey-Wilson polynomials, which themselves are $q$-analogues of the
classical (Hermite, Laguerre, Jacobi) orthogonal polynomials.

Based on an analogy with Macdonald polynomials associated to general root
systems, Macdonald made three conjectures for the Koornwinder polynomials.
In addition to conjectured formulas for principal specialization
\[
k^{(n)}_\lambda(t_0{:}t_1,t_2,t_3;q,t):=
K^{(n)}_\lambda(\dots t^{n-i}t_0\dots;t_0,t_1,t_2,t_3;q,t)
\]
and for the norm with respect to the above inner product, Macdonald made a
third conjecture, which we will call evaluation symmetry, stating that
\beq
\frac{K^{(n)}_\lambda(\dots q^{\mu_i} t^{n-i}t_0\dots;t_0,t_1,t_2,t_3;q,t)}
     {K^{(n)}_\lambda(\dots t^{n-i}t_0\dots;t_0,t_1,t_2,t_3;q,t)}
=
\frac{K^{(n)}_\mu(\dots q^{\lambda_i} t^{n-i}\hat{t}_0\dots;\hat{t}_0,\hat{t}_1,\hat{t}_2,\hat{t}_3;q,t)}
     {K^{(n)}_\mu(\dots
       t^{n-i}\hat{t}_0\dots;\hat{t}_0,\hat{t}_1,\hat{t}_2,\hat{t}_3;q,t)},
\eeq
for suitably modified parameters $\hat{t}_i$.  In \cite{vanDiejenJF:1996},
van Diejen showed that these conjectures were equivalent; evaluation
symmetry was then proved by Sahi \cite{SahiS:1999}, extending work of
Cherednik \cite{CherednikI:1995} for other root systems, using the relevant
``double affine Hecke algebra'' \cite{NoumiM:1994}; see for instance the
book \cite{MacdonaldIG:2003} (which treats all three conjectures directly
via the double affine Hecke algebra).  Essentially, this approach involves
a certain large family (the affine Hecke algebra) of $q$-difference
operators for which the Koornwinder polynomials are eigenfunctions; it also
constructs an associated family of non-symmetric orthogonal polynomials.
(A different approach, also non-symmetric and applicable to arbitrary root
systems, was recently developed by Chalykh \cite{ChalykhOA:2002}.)

In recent work \cite{bcpoly}, we developed a radically different approach
to understanding Koornwinder polynomials (and in particular proving
Macdonald's conjectures).  This approach is in many respects weaker--at
present, it cannot handle the non-symmetric Koornwinder polynomials, and
only works for the root system $BC_n$ (the hardest case for the other
approaches!)--but has a significant advantage in one important respect: it
can be generalized (fairly) easily to the {\em elliptic} level
\cite{bctheta}.  (See also the contributions by Gustafson and Spiridonov to
this volume for discussions of related elliptic special functions.)  This
approach is based on Okounkov's interpolation polynomials
\cite{OkounkovA:1998}, as well as a certain $q$-difference operator that
acts nicely on these polynomials and the Koornwinder polynomials; note this
operator is not, in fact, an element of the affine Hecke algebra, although
it can presumably be constructed using the related theory of raising
operators \cite{MacdonaldIG:2003}.

In \cite{xforms}, inspired by Okounkov's use of an integral operator to
study and construct interpolation polynomials, we gave an explicit
construction of the elliptic analogue of Koornwinder polynomials, using a
sequence of difference and integral operators.  There is thus a
corresponding construction of Koornwinder polynomials obtained by
degenerating from the elliptic case.  In the present note, we describe this
construction, and use it to give yet another proof of two of the three
Macdonald conjectures (principal specialization and norm).

\subsection*{Acknowledgements}
This paper is based on a talk the author gave at the Workshop on Jack,
Hall-Littlewood and Macdonald Polynomials held at the International Centre
for Mathematical Sciences, September 23 through 26, 2003.  The author would
like to thank the organizers for inviting him to that stimulating meeting,
as well as the other participants for making the meeting stimulating.

\subsection*{Notation}  Following \cite{bcpoly}, we define three multivariate
analogues of $q$-symbols:
{\allowdisplaybreaks
\begin{align}
C^0_{\lambda}(x;q,t)
&:=
\prod_{1\le i}(t^{1-i} x;q)_{\lambda_i}
\\
C^-_{\lambda}(x;q,t)
&:=
\prod_{1\le i\le j}
\frac{(t^{j-i} x;q)_{\lambda_i-\lambda_{j+1}}}
     {(t^{j-i} x;q)_{\lambda_i-\lambda_j}}
\\
C^+_{\lambda}(x;q,t)
&:=
\prod_{1\le i\le j}
\frac{(t^{2-i-j} x;q)_{\lambda_i+\lambda_j}}
     {(t^{2-i-j} x;q)_{\lambda_i+\lambda_{j+1}}},
\end{align}
}
with the usual conventions representing products of $C$ symbols via
multiple arguments.  We refer the reader to \cite[\S 2]{bcpoly} for further
discussion of these symbols and the transformations they satisfy.  We also
follow \cite{xforms} in defining two particularly important combinations of
$C$ symbols:
\begin{align}
\Delta^0_\lambda(a|b_1,\dots,b_{2m};q,t)
&=
\left(\frac{(qa)^m}{b_1\cdots b_{2m}}\right)^{|\lambda|}
\frac{C^0_\lambda(b_1,\dots,b_{2m};q,t)}
     {C^0_\lambda(qa/b_1,\dots,qa/b_{2m};q,t)}\\
\Delta_\lambda(a|b_1,\dots,b_{2m};q,t)
&=
\Delta^0_\lambda(a|b_1,\dots,b_{2m};q,t)
\frac{t^{2n(\lambda)} (t/qa)^{|\lambda|} C^0_{2\lambda^2}(qa;q,t)}
     {C^-_\lambda(q,t;q,t) C^+_\lambda(a,qa/t;q,t)}.
\end{align}
(These are, of course, limits of the corresponding symbols of
\cite{xforms} appropriate to the Koornwinder degeneration.)

Given a partition $\lambda$ with at most $n$ parts, the $BC_n$-symmetric
monomial function $m_\lambda$ is defined to be the symmetrization of the
monomial $\prod_i z_i^{\lambda_i}$.  Note that in terms of the usual
monomial symmetric function, we have
$
m_\lambda(\dots z_i\dots) = m_\lambda(\dots z_i+z_i^{-1}\dots).
$
We define a $BC_n$-symmetric function $e_\lambda$ analogously, for
partitions with $\lambda_1\le n$.

If $f$ is a $BC_n$-symmetric polynomial, we define
\begin{align}
\la f\ra^{\prime(n)}_{t_0,t_1,t_2,t_3;q,t}
=
\frac{(q;q)^n}{(t;q)^n 2^n n!}
\int
f(\dots z_i\dots)&
\Delta^{(n)}(\dots z_i\dots;t_0,t_1,t_2,t_3;q,t)\\
&\prod_{1\le i\le n} \frac{dz_i}{2\pi\sqrt{-1}z_i},\notag
\end{align}
and
\[
\la f\ra^{(n)}_{t_0,t_1,t_2,t_3;q,t}
=
\frac{\la f\ra^{\prime(n)}_{t_0,t_1,t_2,t_3;q,t}}
     {\la 1\ra^{\prime(n)}_{t_0,t_1,t_2,t_3;q,t}},
\]
surpressing $(n)$ when it follows from context.  If
$|t_0|,|t_1|,|t_2|,|t_3|,|q|,|t|<1$, then the contour of integration will
be the unit torus; otherwise, the contour needs to be modified to
meromorphically continue from this case.  We finally define
\[
N^{(n)}_\lambda(t_0{:}t_1,t_2,t_3;q,t)
:=
\frac{
\la 
K^{(n)}_\lambda(;t_0,t_1,t_2,t_3;q,t)
K^{(n)}_\lambda(;t_0,t_1,t_2,t_3;q,t)
\ra_{t_0,t_1,t_2,t_3;q,t}
}
{
k^{(n)}_\lambda(t_0{:}t_1,t_2,t_3;q,t)
k^{(n)}_\lambda(t_0{:}t_1,t_2,t_3;q,t)
}
\]

\section{Difference operators}

Of course, the first thing to consider when studying a nice family of
orthogonal polynomials is the normalization of the inner product density
itself.  In the case of the Koornwinder polynomials, this normalization was
given by the following theorem of Gustafson.

\begin{thm}\label{thm:normalization} \cite{GustafsonRA:1990}
For arbitrary complex parameters $q$, $t$, $t_0$, $t_1$, $t_2$, $t_3$ all of
absolute value less than 1,
\[
\la 1\ra^{\prime(n)}_{t_0,t_1,t_2,t_3;q,t}
=
\prod_{1\le j\le n}
\frac{(t^{2n-j-1} t_0t_1t_2t_3;q)}{(t^{n-j+1};q)
\prod_{0\le r<s\le 3} (t^{n-j}t_rt_s;q)
}.
\]
\end{thm}

We will discuss Gustafson's proof in the sequel, but for the present the
following proof will be more relevant.  First, a lemma.

\begin{lem}\label{lem:normalization}
For arbitrary complex numbers $t$, $t_0$, $t_1$,
\[
\prod_{1\le i\le n} (1+R(z_i))
\prod_{1\le i\le n} \frac{(1-t_0z_i)(1-t_1z_i)}{1-z_i^2}
\prod_{1\le i<j\le n} \frac{1-tz_iz_j}{1-z_iz_j}
=
\prod_{1\le i\le n} (1-t^{n-i}t_0t_1),
\]
where $R(z)$ is the operator defined by $R(z)f(z)=f(1/z)$.
\end{lem}

\begin{proof}
If we multiply the left-hand side by the fully antisymmetric polynomial
\[
\Delta(z)=
\prod_{1\le i\le n} (z_i-1/z_i)
\prod_{1\le i<j\le n} (z_i+1/z_i-z_j-1/z_j),
\]
we obtain the sum
\[
\prod_{1\le i\le n}   (1-R(z_i))
\prod_{1\le i\le n}   \frac{-(1-t_0z_i)(1-t_1z_i)}{z_i}
\prod_{1\le i<j\le n} \frac{(1-tz_iz_j)(z_j-z_i)}{z_iz_j}.
\]
Since each term is a Laurent polynomial, the sum is itself a Laurent
polynomial.  Moreover, since
\[
\prod_{1\le i\le n}   \frac{-(1-t_0z_i)(1-t_1z_i)}{z_i}
\prod_{1\le i<j\le n} \frac{(1-tz_iz_j)(z_j-z_i)}{z_iz_j}
\]
is antisymmetric under permutations of the variables, and the group
generated by the $R(z_i)$ is normalized by $S_n$, it follows that the sum
will also be antisymmetric.  Since it is also antisymmetric under each
$R(z_i)$, we find that it is antisymmetric under the full action of
$BC_n$.  But then it must be a multiple of $\Delta(z)$.  Comparing degrees,
we find that the original left-hand side sums to a constant.

To compute this constant, we can proceed in either of two ways.  First, if
we specialize $z_i=t^{n-i}t_0$, only one term on the left survives, which
immediately simplifies to give the desired result.  Alternatively, we can
simply compute the coefficient of the leading monomial of
\[
\prod_{1\le i\le n}   (1-R(z_i))
\prod_{1\le i\le n}   \frac{-(1-t_0z_i)(1-t_1z_i)}{z_i}
\prod_{1\le i<j\le n} \frac{(1-tz_iz_j)(z_j-z_i)}{z_iz_j}.
\]
\end{proof}

We can now give the associated proof of Theorem \ref{thm:normalization}.

\begin{proof}
Factor the integrand as
\begin{align}
\Delta^{(n)}(\dots z_i\dots;t_0,t_1,t_2,t_3;q,t)
&=
\Delta^{(n)}_+(\dots z_i\dots;t_0,t_1,t_2,t_3;q,t)\\
&\phantom{{}={}}
\Delta^{(n)}_+(\dots z_i^{-1}\dots;t_0,t_1,t_2,t_3;q,t),\notag
\end{align}
where
\[
\Delta^{(n)}_+(\dots z_i\dots;t_0,t_1,t_2,t_3;q,t)
=
\prod_{1\le i\le n} \frac{(z_i^2;q)}
{(t_0 z_i,t_1 z_i,t_2 z_i,t_3 z_i;q)}
\prod_{1\le i<j\le n} \frac{(z_i z_j^{\pm 1};q)}{(t z_i z_j^{\pm 1};q)},
\]
and consider the integral
\begin{align}
\int
&\Delta^{(n)}_+(\dots q^{1/2} z_i\dots;t'_0,t'_1,t'_2,t'_3;q,t)\\
&\Delta^{(n)}_+(\dots z_i^{-1}\dots;t_0,t_1,t_2,t_3;q,t)
\prod_{1\le i\le n} \frac{dz_i}{2\pi\sqrt{-1}z_i},\notag
\end{align}
where
\[
(t'_0,t'_1,t'_2,t'_3)=(q^{1/2}t_0,q^{1/2}t_1,q^{-1/2}t_2,q^{-1/2}t_3).
\]
Now, Lemma \ref{lem:normalization} can be expressed in the equivalent form
\[
\prod_{1\le i\le n} (1+R(z_i))
\frac{\Delta^{(n)}_+(\dots q^{1/2} z_i\dots;t'_0,t'_1,t'_2,t'_3;q,t)}
     {\Delta^{(n)}_+(\dots z_i\dots;t_0,t_1,t_2,t_3;q,t)}
=
\prod_{1\le i\le n} (1-t^{n-i}t_0t_1),
\]
and we thus conclude that
\begin{align}
\int&
\Delta^{(n)}_+(\dots q^{1/2} z_i\dots;t'_0,t'_1,t'_2,t'_3;q,t)\\
&\Delta^{(n)}_+(\dots z_i^{-1}\dots;t_0,t_1,t_2,t_3;q,t)
\prod_{1\le i\le n} \frac{dz_i}{2\pi\sqrt{-1}z_i}\notag\\
&\phantom{\Delta^{(n)}_+(\dots z_i^{-1}\dots;t_0,t_1)}=
\frac{(t;q)^n n!}{(q;q)^n}
\prod_{1\le i\le n} (1-t^{n-i}t_0t_1)
\la 1\ra^{\prime(n)}_{t_0,t_1,t_2,t_3;q,t}\notag
\end{align}
If we apply the change of variables $z_i\mapsto q^{-1/2}/z_i$, we obtain a
similar simplification; we thus conclude
\[
\la 1\ra^{\prime(n)}_{t_0,t_1,t_2,t_3;q,t}
=
\prod_{1\le i\le n} \left(\frac{1-t^{n-i}t_2t_3/q}{1-t^{n-i}t_0t_1}\right)
\la 1\ra^{\prime(n)}_{q^{1/2}t_0,q^{1/2}t_1,q^{-1/2}t_2,q^{-1/2}t_3;q,t}
\]
Since the desired right-hand side satisfies the same recurrence, and both
sides are invariant under permutations of $t_0$ through $t_3$, we conclude
that the ratio of the two sides of the desired identity is a function only
of $t_0t_1t_2t_3$, $q$, and $t$.

We can then compute this ratio by expanding the limiting case
$t^{n-1}t_2t_3=1$ via residue calculus.
\end{proof}

The key observation is that this proof can be viewed as being based on
adjointness of difference operators.  We define three $q$-difference
operators as follows.

\begin{defn}
Let $t_0$, $t_1$, $t_2$, $t_3$, $q$, $t$ be arbitrary parameters, and
define difference operators acting on $BC_n$-symmetric polynomials as
follows:
\begin{align}
(D^{-(n)}_q(t)f)(\dots z_i\dots)
&=
\prod_{1\le i\le n} (1+R(z_i))
\prod_{1\le i\le n} \frac{z_i}{1-z_i^2}\\
&\phantom{{}={}}
\prod_{1\le i<j\le n} \frac{1-tz_iz_j}{1-z_iz_j}
f(\dots \sqrt{q} z_i\dots)\notag
\\
(D^{(n)}_q(t_0,t_1;t)f)(\dots z_i\dots)
&=
\prod_{1\le i\le n} (1+R(z_i))
\prod_{1\le i\le n} \frac{(1-t_0z_i)(1-t_1z_i)}{1-z_i^2}\\
&\phantom{{}={}}\prod_{1\le i<j\le n} \frac{1-tz_iz_j}{1-z_iz_j}
f(\dots \sqrt{q} z_i\dots)\notag
\end{align}
\begin{align}
(D^{+(n)}_q&(t_0,t_1,t_2,t_3;t)f)(\dots z_i\dots)\\
&=
\prod_{1\le i\le n} (1+R(z_i))
\prod_{1\le i\le n} \frac{(1-t_0z_i)(1-t_1z_i)(1-t_2z_i)(1-t_3z_i)}{z_i(1-z_i^2)}\notag\\
&\phantom{{}={}}\prod_{1\le i<j\le n} \frac{1-tz_iz_j}{1-z_iz_j}
f(\dots \sqrt{q} z_i\dots)\notag
\end{align}
\end{defn}

\begin{thm}
The above difference operators take $BC_n$-symmetric polynomials to
$BC_n$-symmetric polynomials, acting triangularly with respect to dominance
of monomials:
\begin{align}
D^{-(n)}_q(t)m_\lambda
&=
q^{-|\lambda|/2}
\prod_{1\le i\le n} (1-q^{\lambda_i} t^{n-i})
m_{\lambda-1^n}+\text{dominated terms}\\
D^{(n)}_q(t_0,t_1;t)m_\lambda
&=
q^{-|\lambda|/2}
\prod_{1\le i\le n} (1-q^{\lambda_i} t^{n-i} t_0t_1)
m_{\lambda}+\text{dominated terms}\\
D^{+(n)}_q(t_0,t_1,t_2,t_3;t)m_\lambda
&=
q^{-|\lambda|/2} \prod_{1\le i\le n} (1-q^{\lambda_i} t^{n-i} t_0t_1t_2t_3)
m_{\lambda+1^n}\\
&\phantom{{}={}}+\text{dominated terms}.\notag
\end{align}
Furthermore, they satisfy the following adjointness relations with respect
to the Koornwinder inner product:
\[
\la f D^{(n)}_q(t_0,t_1;t) g\ra^{\prime (n)}_{t_0,t_1,t_2,t_3;q,t}
=
\la g D^{(n)}_q(t'_2,t'_3;t) f\ra^{\prime (n)}_{t'_0,t'_1,t'_2,t'_3;q,t},
\]
where
$(t'_0,t'_1,t'_2,t'_3)=(q^{1/2}t_0,q^{1/2}t_1,q^{-1/2}t_2,q^{-1/2}t_3)$, and
\[
\la f D^{+(n)}_q(t_0,t_1,t_2,t_3;t) g\ra^{\prime (n)}_{t_0,t_1,t_2,t_3;q,t}
=
q^{n/2}
\la g D^{-(n)}_q(t) f\ra^{\prime (n)}_{q^{1/2} t_0,q^{1/2} t_1,q^{1/2} t_2,q^{1/2} t_3;q,t}
\]
\end{thm}

\begin{proof}
The same argument used to prove Theorem \ref{thm:normalization} extends
immediately to give adjointness.  Similarly, that the operators take
polynomials to polynomials follows as in the proof of Lemma
\ref{lem:normalization}.  For instance, for $D^-$, we find that after
clearing the denominator, we obtain $\prod_i (1-R(z_i))$ applied to a
Laurent polynomial in which every monomial is dominated by
$\prod_i z_i^{\lambda_i+n-i}$.
When we antisymmetrize and divide the denominator back out, we thus find
that every monomial of the result is dominated by $\prod_i
z_i^{\lambda_i-1}$ as required.  Moreover, the coefficient of that monomial
is readily computed as given above.  (See, for instance, \cite[Theorem
  3.1]{bcpoly}.)
\end{proof}

\begin{cor}\label{cor:koorn_action_diff}
The difference operators act on Koornwinder polynomials as follows.
\begin{align}
D^{-(n)}_q(t)&K^{(n)}_\lambda(;q^{-1/2} t_0,q^{-1/2} t_1,q^{-1/2}
t_2,q^{-1/2} t_3;q,t)\\
&=
q^{-|\lambda|/2} \prod_{1\le i\le n} (1-q^{\lambda_i} t^{n-i})
K^{(n)}_{\lambda-1^n}(;t_0,t_1,t_2,t_3;q,t)\notag\\
D^{(n)}_q(t_0,t_1;t)&K^{(n)}_\lambda(;q^{1/2} t_0,q^{1/2} t_1,q^{-1/2}
t_2,q^{-1/2} t_3;q,t)\\
&=
q^{-|\lambda|/2}
\prod_{1\le i\le n} (1-q^{\lambda_i} t^{n-i} t_0t_1)
K^{(n)}_\lambda(;t_0,t_1,t_2,t_3;q,t)\notag\\
D^{+(n)}_q(t_0,t_1,t_2,t_3;t)&K^{(n)}_\lambda(;q^{1/2} t_0,q^{1/2} t_1,q^{1/2}
t_2,q^{1/2} t_3;q,t)\\
&=
q^{-|\lambda|/2} \prod_{1\le i\le n} (1-q^{\lambda_i} t^{n-i}
t_0t_1t_2t_3)
K^{(n)}_{\lambda+1^n}(;t_0,t_1,t_2,t_3;q,t)\notag
\end{align}
\end{cor}

\begin{rem}
Given this action on Koornwinder polynomials, it is natural to
wonder how our difference operators relate to the theory of double
affine Hecke algebras.  The operator $D^{(n)}_q$ certainly has such
an interpretation, as follows.  There is a diagram automorphism of
the root system $BC_n$ which gives rise to an outer automorphism
of its Weyl group; using this in the standard way gives an operator
corresponding to translation by the (miniscule) weight $(\frac{1}{2},\dots
\frac{1}{2})$ that commutes (modulo a parameter shift) with the usual
commutative subalgebra.  Symmetrizing this gives a difference
operator which, by comparing actions on Koornwinder polynomials,
must equal $D^{(n)}_q$.  Most likely, the operators $D^{+(n)}_q$
and $D^{-(n)}_q$ arise similarly, as analogues of ``shift operators''
(see \cite[\S 5.9]{MacdonaldIG:2003} for the usual version).
\end{rem}

In particular, we see that $D^-$ acts as a lowering operator, and $D^+$
acts as a raising operator.  Moreover, it is clear that we can combine
these ``first-order'' operators in eight different ways to obtain
``second-order'' operators for which the Koornwinder polynomials are
actually eigenfunctions.  These second-order operators all lie in the
center of the affine Hecke algebra; the first-order operators do not, but
can presumably still be obtained from that theory.

For our present purposes, the main consequence of this action is the
following recurrences for the principal specialization and the norm:

\begin{cor}\label{cor:recur_diff}
For the principal specialization, we have:
\[
\frac{k^{(n)}_\lambda(t_0{:}t_1,t_2,t_3;q,t)}
     {k^{(n)}_\lambda(q^{1/2} t_0{:}q^{1/2}
       t_1,q^{-1/2} t_2,q^{-1/2} t_3;q,t)}
=
q^{|\lambda|/2}
\prod_{1\le i\le n}
\frac{1-t^{n-i}t_0t_1}
     {1-q^{\lambda_i} t^{n-i}t_0t_1}
\]and\begin{align}
&\frac{k^{(n)}_{\lambda+1^n}(t_0{:}t_1,t_2,t_3;q,t)}
     {k^{(n)}_\lambda(q^{1/2} t_0{:}q^{1/2} t_1,q^{1/2} t_2,q^{1/2}
       t_3;q,t)}
\\
&\qquad\qquad
=
q^{|\lambda|/2}
\prod_{1\le i\le n}
\frac{(1-t^{n-i}t_0t_1)(1-t^{n-i}t_0t_2)(1-t^{n-i}t_0t_3)}
     {t^{n-i}t_0(1-q^{\lambda_i} t^{n-i}t_0t_1t_2t_3)}\notag
\end{align}
For the norms of (normalized) Koornwinder polynomials, we have:
\begin{align}
&\frac{N^{(n)}_\lambda(t_0{:}t_1,t_2,t_3;q,t)}
     {N^{(n)}_\lambda(q^{1/2} t_0{:}q^{1/2} t_1,q^{-1/2} t_2,q^{-1/2}
       t_3;q,t)}
\\
&\qquad\qquad=
\prod_{1\le i\le n}
\frac{q^{-\lambda_i}(1-q^{\lambda_i}t^{n-i}t_0t_1)(1-q^{\lambda_i}t^{n-i}t_2t_3/q)}
{(1-t^{n-i}t_0t_1)(1-t^{n-i}t_2t_3/q)},\notag
\end{align}and\begin{align}
&\frac{N^{(n)}_{(\lambda+1^n)}(t_0{:}t_1,t_2,t_3;q,t)}
     {N^{(n)}_\lambda(q^{1/2}
       t_0{:}q^{1/2}t_1,q^{1/2}t_2,q^{1/2}t_3;q,t)}
\\
\nopagebreak&\qquad\qquad=
t_0^{2n} t^{n(n-1)}
\prod_{1\le i\le n}
\frac{(1-t^{n-i}t_1t_2)(1-t^{n-i}t_1t_3)(1-t^{n-i}t_2t_3)}
     {(1-t^{n-i}t_0t_1)(1-t^{n-i}t_0t_2)(1-t^{n-i}t_0t_3)}\notag\\
&\phantom{\qquad\qquad=t_0^{2n} t^{n(n-1)}}
\prod_{1\le i\le n}
\frac{q^{-\lambda_i}(1-q^{\lambda_i} t^{n-i} q)(1-q^{\lambda_i} t^{n-i}t_0t_1t_2t_3)}
     {(1-t^{2n-i-1}t_0t_1t_2t_3)(1-qt^{2n-i-1}t_0t_1t_2t_3)}.
\notag
\end{align}
\end{cor}

\begin{proof}
For the first two recurrence relations, we observe that
$D^{(n)}_q(t_0,t_1;t)$ and $D^{+(n)}_q(t_0,t_1,t_2,t_3;t)$ respect
principal specialization (relative to $t_0$), and thus these relations
follow immediately from the action of these operators on Koornwinder
polynomials.  Similarly, the norm recurrence follows from this action
together with adjointness.
\end{proof}

These recurrences are not quite enough to completely specify these
quantities; there is still freedom when $\lambda_n=0$ to multiply by an
arbitrary function of $t_0t_1t_2t_3$, $q$, and $t$.  To eliminate this
freedom, we will use another, dual, collection of recurrences.

\section{Integral operators}

Gustafson's original proof of Theorem \ref{thm:normalization} was based on
the following integral identity.

\begin{thm}\label{thm:typeI} \cite{GustafsonRA:1994}
For any integer $n\ge 0$, choose complex parameters $q$, $t_0$,\dots,
$t_{2n+1}$, $|q|<1$, such that the sets
\[
\{ q^k t_r:k\ge 0,0\le r\le 2n+1\}
\text{ and }
\{ q^{-k}/t_r:k\ge 0,0\le r\le 2n+1\}
\]
are disjoint, and thus one can choose a contour $C$ containing the first
set and excluding the second set.  Then
\begin{align}
\frac{(q;q)^n}{2^n n!}
\int_{C^n}
\prod_{1\le i<j\le n} (z_i^{\pm 1}z_j^{\pm 1};q)
&\prod_{1\le i\le n} \frac{(z_i^{\pm 2};q)}{\prod_{0\le r<2n+2} (t_r
  z_i^{\pm 1};q)}\frac{dz_i}{2\pi\sqrt{-1}z_i}\\
&\qquad\qquad\qquad\qquad=
\frac{(t_0t_1\cdots t_{2n+1};q)}{\prod_{0\le r<s\le 2n+1} (t_rt_s;q)}.\notag
\end{align}
\end{thm}

\begin{rem}
In addition to the proof in \cite[\S 7]{GustafsonRA:1994}, based on a
multivariate bilateral hypergeometric summation identity, and a proof along
the lines of \cite{xforms} using the fact that when $n$ pairs of parameters
multiply to $q$, the result is a determinant of Askey-Wilson integrals, we
remark that there is a third proof based on the identity
\[
\prod_{1\le i\le n} (1+R(z_i))
\frac{(1-t_0 z_i)\cdots (1-t_n z_i)}{1-z_i^2}
\prod_{1\le i<j\le n} \frac{1}{1-z_iz_j}
=
\prod_{0\le i<j\le n} (1-t_it_j),
\]
which gives an argument along the lines of our proof of Theorem
\ref{thm:normalization} above.  As in that case, this gives rise to pairs
of adjoint difference operators acting on $BC_n$-symmetric polynomials; it
is not clear, however, what significance these operators might have.
\end{rem}

Gustafson's proof of Theorem \ref{thm:normalization} is based on the
following double integral:
\begin{align}
\int_{C^{n+1}} \int_{C^{\prime n}}
\frac{
  \prod_{0\le i<j\le n} (x_i^{\pm 1}x_j^{\pm 1};q)
  \prod_{1\le i<j\le n} (y_i^{\pm 1}y_j^{\pm 1};q)
}{
\prod_{\substack{0\le i\le n\\1\le j\le n}} (\sqrt{t} x_i^{\pm 1} y_j^{\pm 1};q)}
\prod_{1\le i\le n} (y_i^{\pm 2};q) \frac{dy_i}{2\pi\sqrt{-1}y_i}&\\
\prod_{0\le i\le n}
\frac{(x_i^{\pm 2};q)}
     {\prod_{0\le r\le 3} (t_r x_i^{\pm 1};q)}
\frac{dx_i}{2\pi\sqrt{-1}x_i}&,\notag
\end{align}
with appropriate choices of contour.  Both the $x$ and $y$ variables
independently can be integrated out via Theorem \ref{thm:typeI}; the
resulting identity gives a recurrence in $n$ for the Koornwinder
normalization, from which Theorem \ref{thm:normalization} follows
immediately.

Just as the first proof above gives rise to adjoint pairs of difference
operators, Gustafson's proof gives rise to adjoint pairs of {\em integral}
operators.  Defining the operators and proving adjointness is
straightforward; the main difficulty is simply proving that they take
$BC_n$-symmetric polynomials to $BC_n$-symmetric polynomials.  The key fact
is the following generalization of Theorem \ref{thm:typeI}.  Define an
integral operator $I^{*(n)}(q)$ taking $BC_n$-symmetric polynomials to
$S_{2n+2}$-symmetric functions by
\[
(I^{*(n)}(q)f)(t_0,t_1,\dots, t_{2n+1})
=
\int_{C^n} f(\dots z_i\dots) \kappa(\dots z_i\dots) \frac{dz_i}{2\pi
  \sqrt{-1} z_i},
\]
where
\begin{align}
\kappa(\dots z_i\dots) 
={}&
\frac{\prod_{0\le r<s\le 2n+1} (t_rt_s;q)}{(t_0t_1\cdots t_{2n+1};q)}
\frac{(q;q)^n}{2^n n!}\\
&\prod_{1\le i<j\le n} (z_i^{\pm 1}z_j^{\pm 1};q)
\prod_{1\le i\le n} \frac{(z_i^{\pm 2};q)}{\prod_{0\le r<2n+2} (t_r
  z_i^{\pm 1};q)}\notag
\end{align}
and the contour $C$ is as above.

\begin{thm}\label{thm:AWtrans}
If
\[
f(\dots z_i\dots)
=
\prod_{\substack{1\le i\le n\\1\le j\le m}} (y_j+y_j^{-1}-z_i-z_i^{-1}),
\]
then
\begin{align}
(I^{*(n)}(q)f)(t_0,t_1,\dots,t_{2n+1})
={}&
(t_0t_1\cdots t_{2n+1};q)_m^{-1}\\
&\prod_{1\le i\le m} (1+R(y_i))
\frac{\prod_{0\le r<2n+2} (1-t_ry_i)}
     {y_i^n(1-y_i^2)}\notag\\
&\prod_{1\le i<j\le m} \frac{1-q y_iy_j}{1-y_iy_j}
.\notag
\end{align}
\end{thm}

\begin{proof}
The key step is the following lemma.

\begin{lem}
For any $BC_n$-symmetric polynomial $f$,
\begin{align}
(1-R(y))
y^{-n} &\prod_{0\le r<2n}(1-t_ry)
\ (I^{*(n)}(q)f)(t_0,t_1,\dots,t_{2n-1},qy,1/y)\notag\\
&=
(1-t_0t_1\dots t_{2n-1})
y^{-1}(1-y^2)\ 
(I^{*(n-1)}(q)\tilde{f})(t_0,t_1,\dots,t_{2n-1}),
\end{align}
where
\[
\tilde{f}(z_1,\dots,z_{n-1})
=
\prod_{1\le i<n} (y+1/y-z_i-1/z_i)
f(z_1,\dots,z_{n-1},y)
\]
\end{lem}

\begin{proof}
In fact, the two integrals on the left have exactly the same integrand, and
thus their difference is controlled entirely by the difference in contours.
This difference is simply that one contour contains $y$ and excludes $1/y$,
while the other contains $1/y$ and excludes $y$.  We can thus expand the
left-hand side via residue calculus; the result follows.
\end{proof}

In particular, the case $(n,m)$ of the theorem implies the case
$(n-1,m+1)$; since the case $m=0$ is just Theorem \ref{thm:typeI}, the
result follows.
\end{proof}

Note that aside from the factor
$
(t_0t_1\cdots t_{2n+1};q)_m^{-1},
$
the right-hand side is polynomial in $t_0$,\dots, $t_{2n+1}$, and thus the
following three integral operators take $BC_n$-symmetric polynomials to
$BC_{n'}$ symmetric polynomials, for $n'=n+1$, $n$, or $n-1$ as appropriate.

\begin{defn}
Define three integral operators acting on $BC_n$-symmetric polynomials as
follows.
\begin{align}
(I^{+(n)}_t(q)f)(z_1,\dots,z_{n+1})
&=
(I^{*(n)}_t(q)f)(\dots\sqrt{t} z_i^{\pm 1}\dots)\\
(I^{(n)}_t(t_0,t_1;q)f)(z_1,\dots,z_n)
&=
(I^{*(n)}_t(q)f)(t_0,t_1,\dots\sqrt{t} z_i^{\pm 1}\dots)\\
(I^{-(n)}_t(t_0,t_1,t_2,t_3;q)f)(z_1,\dots,z_{n-1})
&=
(I^{*(n)}_t(q)f)(t_0,t_1,t_2,t_3,\dots\sqrt{t} z_i^{\pm 1}\dots).
\end{align}
\end{defn}

\begin{thm}\label{thm:intop_act}
The above operators act on ($BC_n$-symmetric) monomials as follows.
\begin{align}
I^{+(n)}_t(q)m^{(n)}_\lambda
&=
t^{|\lambda|/2}
\prod_{1\le i\le m} 
\frac{1-t^{n+1-\lambda'_i}q^{i-1}}
     {1-t^{n+1}q^{i-1}}
m^{(n+1)}_\lambda + \text{dominated terms},\\
I^{(n)}_t(t_0,t_1;q)m^{(n)}_\lambda
&=
t^{|\lambda|/2}
\prod_{1\le i\le m} 
\frac{1-t^{n-\lambda'_i}q^{i-1}t_0t_1}
     {1-t^{n}q^{i-1}t_0t_1}
m^{(n)}_\lambda + \text{dominated terms},
\end{align}
\begin{align}
I^{-(n)}_t(t_0,t_1,t_2,t_3;q)m^{(n)}_\lambda
&=
t^{|\lambda|/2}
\prod_{1\le i\le m} 
\frac{1-t^{n-1-\lambda'_i}q^{i-1}t_0t_1t_2t_3}
     {1-t^{n-1}q^{i-1}t_0t_1t_2t_3}
m^{(n-1)}_\lambda\\
&\phantom{{}={}} + \text{dominated terms},\notag
\end{align}
where $m^{(n)}_\lambda:=0$ if $\ell(\lambda)>n$.
\end{thm}

\begin{proof}
We first observe that our integral operators are in fact very closely
related to our difference operators; an integral operator acting on the
$z$ variables of a product
\[
f_{n,m}=\prod_{1\le i\le n,1\le j\le m} (y_j+1/y_j-z_i-1/z_i)
\]
becomes a difference operator acting on the $y$ variables of a
corresponding product.  More precisely, we have the following special cases
of Theorem \ref{thm:AWtrans}.
\begin{align}
I^{+(n)}_t(q)_z f_{n,m}
&=
(t^{n+1};q)_m^{-1}
t^{m(n+1)/2}
D^{-(m)}_t(q)_y f_{n+1,m}
\\
I^{(n)}_t(t_0,t_1;q)_z f_{n,m}
&=
(t^n t_0t_1;q)_m^{-1}
t^{mn/2}
D^{(m)}_t(t_0,t_1;q)_y f_{n,m}
\\
I^{-(n)}_t(t_0,t_1;q)_z f_{n,m} 
&=
(t^{n-1} t_0t_1t_2t_3;q)_m^{-1}
t^{m(n-1)/2}
D^{+(m)}_t(t_0,t_1,t_2,t_3;q)_y
f_{n-1,m}.
\end{align}

Now, the product $f_{n,m}$ behaves nicely with respect to dominance of
monomials: we have an expansion
\[
f_{n,m}
=
\sum_{\lambda\subset m^n}
(-1)^{|\lambda|}
m_\lambda(z_1,\dots,z_n)
m_{n^m-\lambda'}(y_1,\dots,y_m)
+
\text{dominated terms},
\]
in the sense that the coefficient of $m_\lambda(z)$ has dominant monomial
$(-1)^{|\lambda|}m_{n^m-\lambda'}(y)$, and vice versa.  (Note that
$n^m-\lambda'$ dominates $\mu$ if and only if $m^n-\mu'$ dominates
$\lambda$, so this condition is indeed symmetrical.)

Thus the fact that the difference operators are triangular implies that the
integral operators are triangular, and similarly for determining the
diagonal coefficients; the theorem follows.
\end{proof}

\begin{lem}
The integral operators satisfy the adjointness relations
\[
\la
g I^{(n)}_t(t_0,t_1;q) f
\ra^{(n)}_{t'_0,t'_1,t'_2,t'_3;q,t}
=
\la
f I^{(n)}_t(t'_2,t'_3;q) g
\ra^{(n)}_{t_0,t_1,t_2,t_3;q,t}
\]
where
$(t'_0,t'_1,t'_2,t'_3)=(t^{1/2}t_0,t^{1/2}t_1,t^{-1/2}t_2,t^{-1/2}t_3)$,
and
\[
\la
h I^{-(n)}_t(t_0,t_1,t_2,t_3;q) f
\ra^{(n-1)}_{t^{1/2} t_0,t^{1/2} t_1,t^{1/2} t_2,t^{1/2} t_3;q,t}
=
\la
f I^{+(n-1)}_t(q) h
\ra^{(n)}_{t_0,t_1,t_2,t_3;q,t},
\]
for any $BC_n$-symmetric polynomials $f$ and $g$, and any
$BC_{n-1}$-symmetric polynomial $h$.
\end{lem}

\begin{proof}
Simply change order of integration.
\end{proof}

\begin{rem}
Note that here we are using the normalized inner product.
\end{rem}

\begin{cor}
The integral operators act on Koornwinder polynomials as follows.
\begin{align}
I^{+(n)}_t(q)
K^{(n)}_\lambda(&;t_0,t_1,t_2,t_3;q,t)\\
=
t^{|\lambda|/2}&
\prod_{1\le i\le m} 
\frac{1-t^{n+1-\lambda'_i}q^{i-1}}
     {1-t^{n+1}q^{i-1}}\notag\\
&
K^{(n+1)}_\lambda(;t^{-1/2} t_0,t^{-1/2} t_1,t^{-1/2} t_2,t^{-1/2} t_3;q,t)
\notag
\\
I^{(n)}_t(t_0,t_1;q)
K^{(n)}_\lambda(&;t_0,t_1,t_2,t_3;q,t)\\
=
t^{|\lambda|/2}&
\prod_{1\le i\le m} 
\frac{1-t^{n-\lambda'_i}q^{i-1} t_0t_1}
     {1-t^{n}q^{i-1} t_0t_1}\notag\\
&K^{(n)}_\lambda(;t^{1/2} t_0,t^{1/2} t_1,t^{-1/2} t_2,t^{-1/2} t_3;q,t)\notag
\\
I^{-(n)}_t(t_0,t_1,t_2,t_3;q)
K^{(n)}_\lambda(&;t_0,t_1,t_2,t_3;q,t)\\
=
t^{|\lambda|/2}&
\prod_{1\le i\le m} 
\frac{1-t^{n-1-\lambda'_i}q^{i-1} t_0t_1t_2t_3}
     {1-t^{n-1}q^{i-1} t_0t_1t_2t_3}\notag\\
&K^{(n-1)}_\lambda(;t^{1/2} t_0,t^{1/2} t_1,t^{1/2} t_2,t^{1/2} t_3;q,t).\notag
\end{align}
\end{cor}

\begin{rem}
In particular, note that
\begin{align}
D^{-(n+1)}_q(t)I^{+(n)}_t(q)
=
I^{-(n)}_t(t_0,t_1,t_2,t_3;t)D^{+(n)}_q(t_0,t_1,t_2,t_3;t)
=
0.
\end{align}
\end{rem}

\begin{cor}\label{cor:recur_int}
For the principal specialization, we have:
\begin{align}
\frac{k^{(n)}_\lambda(t^{1/2} t_0{:}t^{1/2} t_1,t^{-1/2} t_2,t^{-1/2} t_3;q,t)}
     {k^{(n)}_\lambda(t_0{:}t_1,t_2,t_3;q,t)}
&=
t^{-|\lambda|/2}
\prod_{1\le i\le m} 
\frac{1-t^{n}q^{i-1} t_0t_1}
     {1-t^{n-\lambda'_i}q^{i-1} t_0t_1}
\\
\noalign{\noindent and}
\frac{k^{(n+1)}_\lambda(t^{-1/2} t_0{:}t^{-1/2} t_1,t^{-1/2} t_2,t^{-1/2} t_3;q,t)}
     {k^{(n)}_\lambda(t_0{:}t_1,t_2,t_3;q,t)}
&=
t^{-|\lambda|/2}
\prod_{1\le i\le m} 
\frac{1-t^{n+1}q^{i-1}}
     {1-t^{n+1-\lambda'_i}q^{i-1}}
\end{align}
For the norms of (normalized) Koornwinder polynomials, we have:
\begin{align}
&\frac{N^{(n)}_\lambda(t^{1/2}t_0{:}t^{1/2}t_1,t^{-1/2}t_2,t^{-1/2}t_3;q,t)}
     {N^{(n)}_\lambda(t_0{:}t_1,t_2,t_3;q,t)}\\
&\qquad\qquad=
t^{|\lambda|}
\prod_{1\le i\le m} 
\frac{(1-t^{n-\lambda'_i}q^{i-1} t_0t_1)(1-t^{n-\lambda'_i}q^{i-1} t_2t_3/t)}
     {(1-t^{n}q^{i-1} t_0t_1)(1-t^{n}q^{i-1} t_2t_3/t)}\notag\\
\noalign{\noindent and}
&\frac{N^{(n+1)}_\lambda(t^{-1/2} t_0{:}t^{-1/2} t_1,t^{-1/2} t_2,t^{-1/2} t_3;q,t)}
     {N^{(n)}_\lambda(t_0{:}t_1,t_2,t_3;q,t)}\\
&\qquad\qquad=
t^{|\lambda|}
\prod_{1\le i\le m} 
\frac{(1-t^{n+1-\lambda'_i}q^{i-1})(1-t^{n-\lambda'_i}q^{i-1} t_0t_1t_2t_3/t^2)}
     {(1-t^{n+1}q^{i-1})(1-t^n q^{i-1} t_0t_1t_2t_3/t^2)}
\notag
\end{align}
\end{cor}

\begin{proof}
The recurrences for the principal specialization follow from the
observation that the limiting integral corresponding to the principal
specialization of $I^{+(n)}_t(;q)f$ or $I^{(n)}_t(t_0,t_1;q)f$ resolves
(modulo symmetry) into a single residue, and is thus simply the principal
specialization of $f$ itself.  The result thus follows from the action of
these operators on Koornwinder polynomials.  The recurrences for the norm
follow immediately from adjointness.
\end{proof}

\section{The difference-integral representation}

\begin{thm}\cite{SahiS:1999,vanDiejenJF:1996}
The principal specialization and norm of Koornwinder polynomials are given
by the following formulas.
\begin{align}
k^{(n)}_\lambda(t_0{:}t_1,t_2,t_3;q,t)
&=
(t_0 t^{n-1})^{-|\lambda|}
t^{n(\lambda)}
\frac{C^0_\lambda(t^n,t^{n-1}t_0t_1,t^{n-1}t_0t_2,t^{n-1}t_0t_3;q,t)}
     {C^-_\lambda(t;q,t) C^+_\lambda(t^{2n-2}t_0t_1t_2t_3/q;q,t)}
\\
N^{(n)}_\lambda(t_0{:}t_1,t_2,t_3;q,t)
&=
\Delta_\lambda(t^{2n-2}t_0t_1t_2t_3/q|t^n,t^{n-1}t_0t_1,t^{n-1}t_0t_2,t^{n-1}t_0t_3;q,t)^{-1}
\end{align}
\end{thm}

\begin{proof}
The recurrences of Corollary \ref{cor:recur_diff} allow us to deduce the
formulas for $\lambda+1^n$ from the formula for $\lambda$; it thus suffices
to consider the case $\lambda_n=0$.  But then the recurrences of Corollary
\ref{cor:recur_int} prove this case, given that the theorem holds in $n-1$
dimensions.  Since the theorem holds for $\lambda=0$, it holds in general.
\end{proof}

The structure of the above induction gives rise to the following
construction of Koornwinder polynomials.

\begin{thm}
Construct a family $\hat{K}^{(n)}_\lambda(;t_0,t_1,t_2,t_3;q,t)$ of
$BC_n$-symmetric polynomials, defined for nonnegative integers $n$ and
partitions $\lambda$ with $\ell(\lambda)\le n$, as follows.
\begin{enumerate}
\item{} $\hat{K}^{(0)}_0(;t_0,t_1,t_2,t_3;q,t)=1$
\item{} For $n>0$, $\lambda_n=0$,
\begin{align}
\hat{K}^{(n)}_\lambda(;t_0,t_1,t_2,t_3;q,t)
=
&t^{-|\lambda|/2}
\prod_{1\le i\le m} 
\frac{1-t^{n}q^{i-1}}
     {1-t^{n-\lambda'_i}q^{i-1}}\\
&I^{+(n-1)}_t(q)
\hat{K}^{(n-1)}_\lambda(;t^{1/2} t_0,t^{1/2} t_1,t^{1/2}t_2,t^{1/2}t_3;q,t)
\notag\end{align}
\item{} For $n>0$, $\lambda_n>0$,
\begin{align}
\hat{K}^{(n)}_\lambda(;t_0,t_1,t_2,t_3;q,t)
&=
q^{(|\lambda|-n)/2} \prod_{1\le i\le n} (1-q^{\lambda_i} t^{n-i}
t_0t_1t_2t_3/q)^{-1}\\
&\phantom{{}={}}
D^{+(n)}_q(t_0,t_1,t_2,t_3;t)\\
&\phantom{{}={}}
\hat{K}^{(n)}_{\lambda-1^n}(;q^{1/2} t_0,q^{1/2} t_1,q^{1/2} t_2,q^{1/2} t_3;q,t).\notag
\end{align}
\end{enumerate}
Then the resulting polynomials are simply the Koornwinder polynomials.
\end{thm}

\begin{rems}
Similarly, one can define a family of polynomials by
\[
\bar{P}^{*(n)}_\lambda(;q,t,s)
=
t^{-|\lambda|/2}
\prod_{1\le i\le m} 
\frac{1-t^{n}q^{i-1}}
     {1-t^{n-\lambda'_i}q^{i-1}}
I^{+(n-1)}_t(q)
\bar{P}^{*(n-1)}_\lambda(;q,t,s)
\]
\begin{align}
\bar{P}^{*(n)}_{\lambda}(x_1,x_2\dots x_n;q,t,s)
&=
\prod_{1\le i\le n} (x_i+x_i^{-1}-s-s^{-1})\\
&\phantom{{}={}}
\bar{P}^{*(n)}_{\lambda-1^n}(x_1,x_2,\dots x_n;q,t,sq).\notag
\end{align}
The resulting polynomials are simply (the symmetric versions of) Okounkov's
interpolation polynomials \cite{OkounkovA:1998}; see also \cite{bcpoly}.
Indeed, this differs from Okounkov's integral representation for these
polynomials only in that our integral operator is defined by a contour
integral, rather than a sum.  When the polynomial is specialized at a point
of the form $q^{\mu_i} t^{n-i} s$, our contour integral becomes a sum over
partitions by residue calculus, and agrees in that case with Okounkov's
$q$-integral.  Thus the above construction for Koornwinder polynomials can
be viewed as an analogue of Okounkov's representation; in fact, these are
both special cases of the construction given in \cite{xforms}
\end{rems}

\begin{rems}
Unfortunately, the above machinery does not appear to give rise to a
similar proof of evaluation symmetry; of course, we can always refer to the
arguments of van Diejen \cite{vanDiejenJF:1996} or Okounkov
\cite{OkounkovA:1998} showing that evaluation symmetry follows from the
principal specialization formula.
\end{rems}

Another straightforward consequence of our machinery is the following
result of Mimachi.

\begin{thm}\cite{MimachiK:2001}
For any integers $m,n\ge 0$,
\begin{align}
\prod_{\substack{1\le i\le n\\1\le j\le m}}
(y_j+y_j^{-1}-x_i-x_i^{-1})
=
\sum_{\lambda\subset m^n}
(-1)^{|\lambda|}&
K^{(n)}_\lambda(x_1,\dots x_n;t_0,t_1,t_2,t_3;q,t)\\
&K^{(m)}_{n^m-\lambda'}(y_1,\dots y_m;t_0,t_1,t_2,t_3;t,q)
\notag\end{align}
\end{thm}

\begin{proof}
Clearly the left-hand side admits {\em some} expansion of the form
\[
\sum_{\lambda,\mu\subset m^n}
c_{\lambda\mu}
K^{(n)}_\lambda(x_1,\dots x_n;t_0,t_1,t_2,t_3;q,t)
K^{(m)}_\mu(y_1,\dots y_m;t_0,t_1,t_2,t_3;t,q).
\]
If we apply one of the ``second-order'' difference operators for which
$K^{(n)}_\lambda$ is a basis of eigenfunctions, the proof of Theorem
\ref{thm:intop_act} turns this composition of two difference operators in
the $x$ variables into a composition of two integral operators in the $y$
variables, for which $K^{(m)}_\mu$ is a basis of eigenfunctions.  Comparing
the two eigenvalues, we find that $c_{\lambda\mu}=0$ unless
$\mu=n^m-\lambda'$.  The coefficient then follows by an examination of
dominant terms.
\end{proof}

\bibliographystyle{plain}

\begin{thebibliography}{10}

\bibitem{AskeyR/WilsonJ:1985}
R.~Askey and J.~Wilson.
\newblock {\em Some basic hypergeometric orthogonal polynomials that generalize
  {J}acobi polynomials}.
\newblock Number 319 in Memoirs of the AMS. Amer. Math. Soc., Providence, RI,
  1985.

\bibitem{ChalykhOA:2002}
O.~A. Chalykh.
\newblock Macdonald polynomials and algebraic integrability.
\newblock {\em Adv. Math.}, 166(2):193--259, 2002.

\bibitem{CherednikI:1995}
I.~Cherednik.
\newblock Double affine {H}ecke algebras and {M}acdonald's conjectures.
\newblock {\em Ann. of Math. (2)}, 141(1):191--216, 1995.

\bibitem{GustafsonRA:1990}
R.~A. Gustafson.
\newblock A generalization of {S}elberg's beta integral.
\newblock {\em Bull. Amer. Math. Soc. (N.S.)}, 22(1):97--105, 1990.

\bibitem{GustafsonRA:1994}
R.~A. Gustafson.
\newblock Some {$q$}-beta and {M}ellin-{B}arnes integrals on compact {L}ie
  groups and {L}ie algebras.
\newblock {\em Trans. Amer. Math. Soc.}, 341(1):69--119, 1994.

\bibitem{KoornwinderTH:1992}
T.~H. Koornwinder.
\newblock {A}skey-{W}ilson polynomials for root systems of type {$BC$}.
\newblock In Donald St.~P. Richards, editor, {\em Hypergeometric functions on
  domains of positivity, Jack polynomials, and applications (Tampa, FL, 1991)},
  Contemp. Math. 138, pages 189--204. Amer. Math. Soc., Providence, RI, 1992.

\bibitem{MacdonaldIG:2003}
I.~G. Macdonald.
\newblock {\em Affine {H}ecke algebras and orthogonal polynomials}, volume 157
  of {\em Cambridge Tracts in Mathematics}.
\newblock Cambridge University Press, Cambridge, 2003.

\bibitem{MimachiK:2001}
K.~Mimachi.
\newblock A duality of {M}acdonald-{K}oornwinder polynomials and its
  application to integral representations.
\newblock {\em Duke Math. J.}, 107(2):265--281, 2001.

\bibitem{NoumiM:1994}
M.~Noumi.
\newblock Macdonald-{K}oornwinder polynomials and affine {H}ecke rings.
\newblock {\em S\=uri\-kaiseki\-kenky\=u\-sho K\=oky\=u\-roku}, (919):44--55,
  1995.
\newblock Various aspects of hypergeometric functions (Japanese) (Kyoto, 1994).

\bibitem{OkounkovA:1998}
A.~Okounkov.
\newblock {$BC$}-type interpolation {M}acdonald polynomials and binomial
  formula for {K}oornwinder polynomials.
\newblock {\em Transform. Groups}, 3(2):181--207, 1998.

\bibitem{bcpoly}
E.~M. Rains.
\newblock {$BC_n$}-symmetric polynomials.
\newblock arxiv:math.QA/0112035.

\bibitem{bctheta}
E.~M. Rains.
\newblock {$BC_n$}-symmetric theta functions.
\newblock arxiv:math.CO/0402113.

\bibitem{xforms}
E.~M. Rains.
\newblock Transformations of elliptic hypergeometric integrals.
\newblock arxiv:math.QA/0309252.

\bibitem{SahiS:1999}
S.~Sahi.
\newblock Nonsymmetric {K}oornwinder polynomials and duality.
\newblock {\em Ann. of Math. (2)}, 150(1):267--282, 1999.

\bibitem{vanDiejenJF:1996}
J.~F. van Diejen.
\newblock Self-dual {K}oornwinder-{M}acdonald polynomials.
\newblock {\em Invent. Math.}, 126(2):319--339, 1996.

\end{thebibliography}

\end{document}